\documentclass[11pt]{article}

\usepackage[latin1]{inputenc}
\usepackage{graphicx}
\usepackage{epsfig}

\newtheorem{theorem}{Theorem}

\newfont{\bb}{msbm10 at 11pt}
\def\r{\hbox{\bb R}}

\begin{document}

\title{A comparison result for  radial solutions of the mean curvature  equation}
\author{Rafael L\'opez\footnote{Partially
supported by MEC-FEDER
 grant no. MTM2007-61775.}\\
Departamento de Geometr\'{\i}a y Topolog\'{\i}a\\
Universidad de Granada\\
18071 Granada (Spain)\\
e-mail: {\tt rcamino@ugr.es}}

\date{}

\maketitle

\begin{abstract} We establish two comparison results between the solutions of a class of mean curvature equations and  pieces of arcs of circles that satisfy the same Neumann boundary condition. Finally we present a number of examples where our estimates can be applied, some of them have a physical motivation.
\end{abstract}


\section{Introduction and statements of results}


In this paper we wish to estimate the solutions of equation of mean curvature type
\begin{equation}\label{r1}
\frac{1}{r}\bigg(r\frac{u'(r)}{\sqrt{1+u'(r)^2}}\bigg)'=f(r),\hspace*{1cm}0\leq
r\leq c
\end{equation}
that satisfies a Neumann boundary  condition
\begin{equation}\label{rn}
u'(c)=\tan(\gamma),\hspace*{1cm}\gamma\in [0,\frac{\pi}{2}).
\end{equation}
As usually, by $'$ we denote the derivative $\frac{d}{dr}$. Exactly we compare the solutions with pieces of arcs of circles, which are the solutions of (\ref{r1})-(\ref{rn}) when $f$ is a constant function. Under appropriate assumptions on the function $f$ we are able to show that the graphic of a solution of (\ref{r1})-(\ref{rn}) can be sandwiched between two arcs of circles with the same boundary condition.  We remark that we do not explicitly address the question of the
existence of solutions of Equation (\ref{r1}).

Equation (\ref{r1}) is the expression in radial coordinates of the prescribed mean curvature equation
\begin{equation}\label{r11}
\mbox{div}\bigg(\frac{Du}{\sqrt{1+|Du|^2}}\bigg)=2
H(x),\hspace*{1cm}x\in\Omega\subset\r^2
\end{equation}
where $\Omega$ is an open set in $\r^2$. In such case,  $H$ is the mean curvature of the
non-parametric surface $z=u(x)$ in   Euclidean three-space $\r^3$. Thus a radial solution  $u=u(r)$, $r=|x|$, defines
 a curve $\alpha(r)=(r,u(r))$ in such way that the surface obtained by rotating $\alpha$ with respect to the $z$-axis has mean curvature $f(r)/2$.
Equation (\ref{r11}) may be used to model a number of important problems in
 mechanics. For example, it appears in the context of the
isoperimetric problem of least surface area bounding a given volume. Under the boundary condition (\ref{rn}), our equation appears in capillary theory  as a mathematical model for
 the equilibrium shape of a liquid surface with constant surface tension in a uniform gravity field with prescribed contact angle with the vertical walls. More examples can seen in Section \ref{applications}.

Returning with the problem (\ref{r1})-(\ref{rn}),  we require that $u$ be a classical solution on
$[0,c]$ and that $f$ is sufficiently smooth. Consider the natural boundary conditions
\begin{equation}\label{r2}
u(0)=u_0,\hspace*{1cm}u'(0)=0.
\end{equation}
and that $I=[0,c]$
is the interval where $u$ is defined. In order to state our results,  we take a piece of circle
with the same slope than $u$ at $r=c$ and that coincides with $u$ at
the origin. Exactly, let
$$y(r)=R+u_0-\sqrt{R^2-r^2},\hspace*{1cm}R=c\ \frac{\sqrt{1+u'(c)^2}}{u'(c)}.$$
The graphic of $y$ is a piece of a lower halfcircle with $y(0)=u_0$ and $y'(0)=0$. The choice of
the radius $R$ is such that  $y'(c)=u'(c)$. Thus $y(r)$ is a solution of (\ref{r1})-(\ref{rn})-(\ref{r2}) for $f(r)=2/R$. In such setting, we
compare $u$ with the circle $y$. Throughout this paper, we  suppose that the next assumption holds for the
function $f$:

\emph{ASSUMPTION.} We suppose that the function $f$ satisfies the
following conditions:
\begin{enumerate}
\item $f(0)\geq 0$,
\item $f$ is an increasing function on $r$ and
\item $f''(r)\geq 0$ for $0\leq r\leq c$.
\end{enumerate}
With the above notation, we state our results.

\begin{theorem}\label{t1} Let $u$ be a solution of
(\ref{r1})-(\ref{r2}) defined in the interval $[0,c]$. Then
$$u(r)<y(r),\hspace*{1cm}0<r\leq c.$$
\end{theorem}

For the next result, we descend vertically the circle $y(r)$ until it touches
with the graphic of $u$ at $r=c$. We call $w=w(r)$ the new position
of $y$, that is, $w(r)=y(r)-y(c)+u(c)$.

\begin{theorem}\label{t2}
Let $u$ be a solution of (\ref{r1})-(\ref{r2}) defined in the
interval $[0,c]$. Then
$$w(r)<u(r),\hspace*{1cm}0\leq r<c.$$
\end{theorem}

\begin{figure}[hbtp]
\includegraphics[width=10cm]{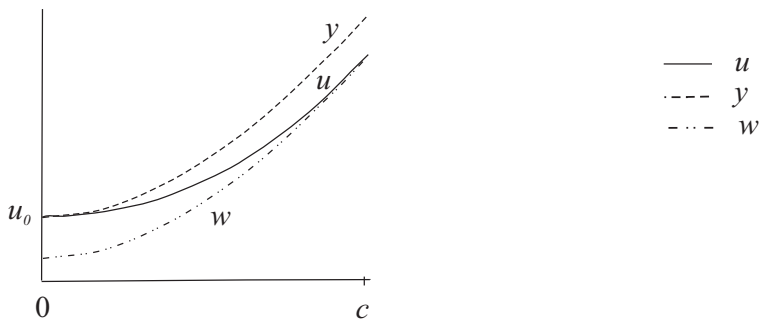}
\caption{The solution $u$ lies sandwiched between the circles $y$ and
$w$.} \label{AML7035-figure1}
\end{figure}

As conclusion, the solution $u$ lies between two pieces of circles,
namely, $y$ and $w$, such that the slopes of the three functions
agree at the points $r=0$ and $r=c$ and the graphic of $u$ coincides with $y$ and $w$ at $r=0$ and $r=c$ respectively. See Figure \ref{AML7035-figure1}. We remark that with appropriate modifications the conclusions of both theorems hold even in the case that the maximal interval of definition of $u$ is  $[0,c)$, where
there exists $\lim_{r\rightarrow c}u(r)$ and $\lim_{r\rightarrow c}\frac{u'(r)}{\sqrt{1+u'(r)^2}}=1$ ($\gamma=\pi/2$).

We point out that the maximum principle for elliptic equations lies behind our proofs.  Actually, we compare the solutions of (\ref{r1})-(\ref{r2}) with the ones of (\ref{r1}) but changing $f(r)$ by a constant $k$. In the latter case, the solution $u$ has the form $\lambda-\sqrt{4/k^2-r^2}$, whose graphic is an arc of circle.  For this reason, our results are within a more general framework with the appellation  of \emph{comparison results} for solutions of divergence structure equations of type $\mbox{div}A(x,u,Du)+B(x,u,Du)=0$ in domains of $\r^2$, where a certain hypothesis of monotonicity of growth is required for $B$. For a wide presentation of the known results in this direction we refer the reader to the classical book of Gilbarg and Trudinger \cite{gt} and an up-to-date modern treatment of the maximum principle of Pucci and Serrin \cite{ps2} (see also \cite{ps1}). For a detailed discussion of comparison principles, see particularly \cite[chapters 2 and 3]{ps2}.

\section{The proofs}
Let $\psi(r)$ be the angle that makes the graphic of $u$ with the
$r$-axis at each point $r$, that is $\tan\psi(r)=u'(r)$. Put
$$\sin\psi(r)=\frac{u'(r)}{\sqrt{1+u'(r)^2}}.$$
Then Equation (\ref{r1}) writes as
\begin{equation}\label{r12}
\frac{1}{r}\bigg(r\sin\psi(r)\bigg)'= f(r).
\end{equation}
A first integration yields
\begin{equation}\label{inte}
\sin\psi(r)=\frac{1}{r}\int_0^r t f(t)\ dt.
\end{equation}
Because the function $f$ is positive, the integrand of (\ref{inte}) is positive too and thus,
$\sin\psi(r)>0$ at $(0,c]$. This means that $u$ is a strictly
increasing function on $r$. Since $f$ is increasing on $r$, fixed a real number $r\in (0,c]$, we have
$f(0)< f(t) < f(r)$ for any $t\in (0,r)$. Putting these inequalities in
the integrand of (\ref{inte}),  we obtain for $0< r\leq c$
\begin{equation}\label{efe}
\frac{r}{2} f(0)< \sin\psi(r)< \frac{r}{2}f(r).
\end{equation}

We work with Equation (\ref{r12}) as follows:
$$\sin\psi(r)+r(\sin\psi(r))'=r f(r)$$
or \begin{equation}\label{kurva}
(\sin\psi(r))'=f(r)-\frac{\sin\psi(r)}{r}. \end{equation}
 The left
side in the above equation, namely, the term $(\sin\psi(r))'$,  is
the curvature $\kappa$ of the generating curve $\alpha$
of the surface $z=u(x)=u(|x|)$, that is,
$$\kappa(r)=(\sin\psi(r))'=\frac{u''(r)}{(1+u'(r)^2)^{3/2}}.$$
The key in our proofs comes from the fact that the function
$\kappa(r)$ is an increasing function on $r$. Exactly, we have
$$\kappa'(r)=(\sin\psi(r))''=f'(r)+2\frac{\sin\psi(r)}{r^2}-\frac{f(r)}{r}.$$
By using the left inequality of (\ref{efe}) and an integration by parts, we conclude
\begin{equation}\label{kurva2}
(\sin\psi)''(r)> f'(r)+\frac{f(0)}{r}-\frac{f(r)}{r}=\frac{1}{r}\int_0 ^r t f''(t)\ dt \geq 0,
\end{equation}
where we have used the fact that $f$ is a convex function.

\emph{Proof of Theorem \ref{t1}.} The angle $\psi^y(r)$ and the curvature
$\kappa^y$ of the function $y(r)$  are respectively
$$\sin\psi^y(r)=\frac{r}{R},\hspace*{1cm}\kappa^y(r)=\frac{1}{R}.$$
We recall that the curvature of the circle $y(r)$ is constant. At $r=0$, we compare the curvatures $\kappa$ and $\kappa^y$. From (\ref{efe})
$$\kappa^y(0)=\frac{1}{R}=\frac{\sin\psi(c)}{c}> \frac{f(0)}{2}.$$
Using this inequality, the expression of $\kappa$ in (\ref{kurva}) and using  the left inequality of (\ref{efe}) again,
we have
$$\kappa(0)=(\sin\psi)'(0)=f(0)-\frac{\sin\psi(r)}{r}(0)\leq \frac{f(0)}{2}<\kappa^y(0).$$
As $\kappa(0)<\kappa^y(0)$, $y(0)=u(0)$ and $y'(0)=u'(0)$,  the graphic of $y$
lies above of $u$ around of $r=0$. Theorem \ref{t1} asserts that
this occurs in the interval $(0,c]$. The proof is by contradiction.
Assume that the graphic of $u$ acroses the graphic of $y$ at some
point. Let $r=\delta\leq c$ the first value where this occurs, that is, $u(r)<y(r)$ for $r\in (0,\delta)$ and $u(\delta)=y(\delta)$.  Then $u'(\delta)\geq y'(\delta)$ and so,
$\sin\psi(\delta)\geq\sin\psi^y(\delta)$. As $u'(0)=y'(0)$, we have
\begin{equation}\label{contr1}
\int_0^\delta\bigg(\kappa(t)-\kappa^y(t)\bigg)\
dt=\int_0^\delta\bigg((\sin\psi(t))'- (\sin\psi^y(t))'\bigg)\
dt=\sin\psi(\delta)-\sin\psi^y(\delta)\geq 0.
\end{equation}
On the other hand, as  $\kappa(0)<\kappa^y(0)$ and  the above integral is non-negative,  there exists $\bar{r}\in
(0,\delta)$ such that $\kappa(\bar{r})>\kappa^y(\bar{r})$.
Because $\kappa$ is increasing on $r$, we have for $r\in [\bar{r},c]$
$$\kappa(r)>\kappa(\bar{r})>\kappa^y(\bar{r})=\kappa^y(r).$$
 Thus and since $\bar{r}\leq\delta\leq c$,
\begin{eqnarray*}
0&<&\int_{\bar{r}}^c \bigg(\kappa(t)-\kappa^y(t)\bigg)\ dt\leq\int_{\delta}^c \bigg(\kappa(t)-\kappa^y (t)\bigg)\ dt\\
&=&\int_{\delta}^c\bigg((\sin\psi(t))'- (\sin\psi^y(t))'\bigg)\
dt=\sin\psi^y(\delta)-\sin\psi(\delta),
\end{eqnarray*}
in contradiction with (\ref{contr1}). This show Theorem \ref{t1}.

\emph{Proof of Theorem \ref{t2}.} Arguing in a similar way,  we compare the curvatures $\kappa^w$ and
$\kappa$ at the point $r=c$. Using the right inequality of (\ref{efe}), we have
$$\kappa(c)=f(c)-\frac{\sin\psi(c)}{c}>\frac{\sin\psi(c)}{c}=\kappa^y(c)=\kappa^w(c).$$
As $\kappa(c)>\kappa^w(c)$, $w(c)=u(c)$ and $w'(c)=u'(c)$, the graphic of $u$ lies above than the circle $w$ around $r=c$. Thus $w(r)<u(r)$ in some  interval $(\delta,c)$.
Again, the proof is by contradiction. We suppose that the
graphic of $w$ acroses the graphic of $u$ at some point. Denote by
$\delta$ the largest number such that $w(r)<u(r)$ for $r\in (\delta,c)$ and $w(\delta)=u(\delta)$.  For this value,
$w'(\delta)=y'(\delta)\leq u'(\delta)$ and $\sin\phi^y(\delta)\leq\sin\psi(\delta)$. Then
\begin{equation}\label{contr2}
\int_{\delta}^c \bigg(\kappa(t)-\kappa^w(t)\bigg) \
dt=\int_{\delta}^c\bigg((\sin\psi(t))'- (\sin\psi^y(t))'\bigg)\
dt=\sin\psi^y(\delta)-\sin\psi(\delta)\leq 0.
\end{equation}
We have used that $u'(c)=w'(c)=y'(c)$. As
$\kappa(c)-\kappa^w(c)>0$ and the integral in (\ref{contr2}) is non-positive, then there would be
$\bar{r}\in(\delta,c)$ such that
$\kappa(\bar{r})<\kappa^w(\bar{r})$. Because $\kappa$ is an
increasing function on $r$, for any $r\in [0,\bar{r}]$
$$\kappa(r)<\kappa(\bar{r})<\kappa^w(\bar{r})=\kappa^w(r).$$
Since $\delta<\bar{r}$,
$$0>\int_0^{\delta} \bigg(\kappa(t)-\kappa^w(t)\bigg)\ dt=\int_0^{\delta}\bigg((\sin\psi(t))'- (\sin\psi^y(t))'\bigg)\
dt=\sin\psi(\delta)-\sin\psi^y(\delta),$$
where we use the fact that $u'(0)=y'(0)$.
 This contradicts the inequality (\ref{contr2}) and proves Theorem \ref{t2}.

\noindent{\bf Remark.} We point out that the assumptions on the function $f$ are not necessary to get our results. Let $f(r)=\sin(r)$. In the interval $(0,\pi)$, this function is positive and increasing on $r$ but $f''<0$. We compute the function $\sin\psi(r)$  for this choice of $f$. From (\ref{inte}),
$$\sin\psi(r)=\frac{\sin(r)}{r}-\cos(r)$$
and the curvature of the graphic of $u$ satisfies
$$\kappa'(r)=(\sin\psi(r))''=\frac{r^2-2}{r^3}(r\cos(r)-\sin(r)).$$
Thus the curvature function $\kappa(r)$ is increasing on $r$ in the interval $(0,\sqrt{2})$ and consequently the statements of Theorems \ref{t1} and \ref{t2} are true for any $c\in (0,\sqrt{2})$. However, in the same interval $(0,\sqrt{2})$, $f''(r)<0$ and thus the inequality in (\ref{kurva2}) yields
$\int_0^r t f''(t)\ dt<0$ for any $r\in (0,\sqrt{2})$. The thing in this case is that the left inequality on (\ref{efe}) is rough too.

\section{Applications}\label{applications}

In this section we apply Theorems \ref{t1} and \ref{t2} into several
specific examples.

\emph{A. Capillary surfaces.} The equation of the equilibrium shape of a liquid surface with constant surface tension in a uniform gravity field is governed by the equation (\ref{r11}) with $2H=B u$. The number  $B$ is a physical constant that is positive or negative according to whether the gravitational field is acting downward or upward. Here we consider $B>0$. In the case that the solution is radial, $f(r)=B u(r)$ and the natural physical boundary condition is  $u'(c)=\tan\gamma$, where $\frac{\pi}{2}-\gamma$ is the angle between the liquid surface and the fixed boundary. Note  that in this situation $f=f(u)$.
Assume that $u(0)=u_0>0$. Since $f(0)=B u_0>0$, we obtain from (\ref{inte}) that $\sin\psi(r)>0$. This proves that $u'(r)>0$ and $f$ is an increasing function on $r$. On the other hand, the sign of $f''(r)$ is the same than $u''(r)$ and $\kappa(r)$. From (\ref{efe}) and  (\ref{kurva}),
$$\kappa(r)=B u(r)-\frac{\sin\psi(r)}{r}>\frac{B u(r)}{2}>\frac{B u_0}{2}>0.$$
As conclusion, we can apply  our results if both $u_0$ and $B$ are positive. This allows to obtain estimates of the volume of the fluid of the liquid drop. This was used in  \cite{fi} to estimate the volume of a capillary surface.

 \emph{B. Capillary for compressible fluids}. In capillarity theory, we take into account the effect of the virtual motions of fluid particles in the internal energy of the fluid. This is caused by the fluid compressibility. Then the equation for the fluid surface height in a capillary circular tube is (\ref{r1}) where the function $f$ is
 $$f(r)=\frac{-a}{\sqrt{1+u'(r)^2}}+b \exp{(au(r))}+c,$$
and  $a>0$, $b>0$ and  $c$ are real numbers ($a$ is called the compressibility constant). Here $f=f(u,u')$. For these values of  $a$ and $b$, one can show that $f$ satisfies our  Assumption. See \cite{fl} for more details. The estimates establish upper and lower bounds for the height
 solutions.

\emph{C. Rotating liquid drops.} In  absence of gravity, we
consider  the steady rigid rotation of a homogeneous incompressible
fluid drop which is surrounded by a rigidly rotating incompressible
fluid. In mechanical equilibrium, we say that the drop is a rotating
liquid drop \cite{bs,we}.  In the case that the drop is asymmetric,
the shape of the interface is locally governed by equation (\ref{r1}) with
$f(r)=a r^2+b$, for constants $a\not=0,b$. The expressions of these constants involve the angular velocity, the density and the surface tension of the fluid of the drop.  For appropriate initial boundary conditions and successive reflections, it is possible to obtain rotating drops homeomorphic to balls, and other ones that adopt toroidal configurations. See  \cite{aq,gu}. On the other hand, this choice of $f$ appears in the study of the motion of a two-fluid interface in a rotating Hele-Shaw  cell. The interface shapes balance the centrifugal and capillary forces \cite{sc,cmco}.

In the case that $a>0 $ and $b\geq 0$, the function $f$ is under the hypothesis of our Assumption.

\emph{D. The function $f$ is linear.} Assume now that $f(r)=ar+b$, where $a$ and $b$ are real numbers. This setting differs from the capillarity theory   where in such case $f$ was a linear function of $u(r)$ (part A of this section). If $a$ and $b$ are non-negative numbers, then we are in conditions of the Assumption on $f$.

\emph{E. The function $f$ is exponential.} Consider $f(r)=a \exp{(r)})$, where $a>0$. Then $f$ lies under the assumptions of our results.

\end{document}